\documentclass[a4paper,12pt]{article}
\usepackage{amsfonts}
\usepackage[xdvi]{graphicx}
\usepackage[dvips]{color}
\usepackage{amsmath,amssymb}

\newcommand{\C}{{\mathbb C}}
\newcommand{\N}{{\mathbb N}}
\newcommand{\Z}{{\mathbb Z}}

\makeatletter

\@addtoreset{equation}{section}

\newcounter{def}[section]
\renewcommand{\thedef}{\stepcounter{def}\thesection.\@arabic\c@def }

\begin{document}
\setlength{\baselineskip}{24pt}
\begin{center}
\textbf{\LARGE{On the degeneration of Kovalevskaya exponents of Laurent series solutions of 
quasi-homogeneous vector fields}}
\end{center}
\vskip.5\baselineskip

\setlength{\baselineskip}{14pt}

\begin{center}
Advanced Institute for Materials Research, Tohoku University,\\
Sendai, 980-8577, Japan
\vskip.5\baselineskip
\large{Hayato CHIBA} \footnote{E mail address : hchiba@tohoku.ac.jp}
\end{center}
\begin{center}
Feb. 1, 2026; \,\, revised June 8, 2026
\end{center}
\vskip.5\baselineskip

\begin{center}
\textbf{Abstract}
\end{center}
A structure of families of Laurent series solutions of a quasi-homogeneous vector field is studied,
where a given vector field is assumed to have a commuting vector field.
For an $m$ dimensional vector field, a family of Laurent series solutions is called principal 
if it contains $m$ arbitrary parameters, and called non-principal if the number is smaller than $m$.
Starting from a principal Laurent series solutions,
a systematic method to obtain a non-principal Laurent series solutions is given.
In particular, from the Kovalevskaya exponents of the principal Laurent series solutions, 
which is one of the invariants of quasi-homogeneous vector fields, 
the Kovalevskaya exponents of the non-principal Laurent series solutions are obtained
by using the commuting vector field.


\section{Introduction}

A differential equation defined on a complex region is said to have the Painlev\'{e} property
if any movable singularity 
(a singularity of a solution which depends on an initial condition) of any solution is a pole.
Among them, an important class is the equations that all solutions are meromorphic.

Painlev\'{e} and his group classified second order ODEs having the Painlev\'{e} property
and found new six differential equations called the Painlev\'{e} equations
$\text{P}_\text{I}, \cdots , \text{P}_\text{VI}$.
Nowadays, it is known that they are written in Hamiltonian forms
\begin{eqnarray*}
(\text{P}_\text{J}): \frac{dq}{dz} = \frac{\partial H_J}{\partial p}, \quad
\frac{dp}{dz} = -\frac{\partial H_J}{\partial q}, \quad J = \text{I}, \cdots, \text{VI}.
\end{eqnarray*}
Among six Painlev\'{e} equations, the Hamiltonian functions of the first, second and fourth
Painlev\'{e} equations are polynomials in both of the independent variable $z$
and the dependent variables $(q,p)$.
They are given by
\begin{eqnarray}
H_{\text{I}}(q,p) &=& \frac{1}{2}p^2 - 2q^3 - zq, \label{1-1}\\
H_\text{II}(q,p) &=& \frac{1}{2}p^2 - \frac{1}{2}q^4 - \frac{1}{2}zq^2 - \alpha q, \nonumber \\
H_\text{IV}(q,p) &=& -pq^2 + p^2q - 2pqz - \alpha p + \beta q, \nonumber
\end{eqnarray}
respectively, where $\alpha, \beta \in \C$ are arbitrary parameters.
Since they satisfy Painlev\'{e} property and the right hand sides are polynomials,
any solutions are meromorphic. 

In general, a polynomial $H(x_1, \cdots ,x_m)$ is called a quasi-homogeneous polynomial if there 
are a tuple of positive integers $(a_1, \cdots ,a_m)$, called the weight, and $h$, called the 
weighted degree, such that
\begin{eqnarray}
H(\lambda ^{a_1}x_1 , \cdots ,\lambda ^{a_n}x_m) = \lambda ^h H(x_1, \cdots ,x_m)
\label{quasi}
\end{eqnarray}
for any $\lambda \in \C$.
The above Hamiltonians $H_{\text{I}}, H_{\text{II}}, H_{\text{IV}}$ are quasi-homogeneous with respect to the weights given in Table 1,
if we ignore terms including parameters $\alpha , \beta$.
In Chiba \cite{Chi4}, possible weights arising from Hamiltonian systems are classified 
from a view point of singularity theory and it is shown that 
they are related to Painlev\'{e} equations.

\begin{table}[h]
\begin{center}
\begin{tabular}{|c||c|c|c|}
\hline
 & $\mathrm{deg}(q,p,z)$ & $\mathrm{deg} (H)$ & $\kappa$ \\ \hline \hline
$\text{P}_\text{I}$  & $(2,3,4)$ & 6 & 6  \\ \hline
$\text{P}_\text{II}$  & $(1,2,2)$ & 4 & 4  \\ \hline
$\text{P}_\text{IV}$  & $(1,1,1)$ & 3 & 3  \\ \hline
\end{tabular}
\end{center}
\caption{ $\mathrm{deg} (H)$ denotes the weighted degree of the Hamiltonian function with
respect to the weight $\mathrm{deg}(q,p,z)$.
$\kappa$ denotes the Kovalevskaya exponent defined in Section 2.}
\end{table}

For a general vector field $F=(f_1, \cdots ,f_m)$ on $\C^m$,
its weight $(a_1, \cdots ,a_m)$ is defined in a similar manner as
\begin{equation}
f_i(\lambda ^{a_1}x_1 ,\cdots , \lambda ^{a_m}x_m) = \lambda ^{a_i+\gamma }f_i(x_1, \cdots ,x_m),
\quad i=1,\cdots ,m, \quad \gamma \in \N,
\end{equation}
if it exists.
Once a weight is given with $\gamma =1$, to construct a Laurent series solution of the equation
$dx_i/dz = f_i$ is straightforward.
It is expressed in the form
\begin{align*}
x_i(z) = c_i(z-\alpha _0)^{-a_i} + \sum^\infty_{j=1} d_{i,j} (z-\alpha _0)^{-a_i+j},
\end{align*}
where coefficients $c_i$ and $d_{i,j}$ are determined by substituting this expression into the equation.
Important features are that a position of a pole $\alpha _0$ can take an arbitrary value
depending on the initial condition, that is called the movable singularity,
and that the order $a_i$ of the pole is the same as the weight of $x_i$.

From Laurent series solutions, we define the Kovalevskaya exponents as follows:
As an example, we consider the first Painlev\'{e} equation (\ref{1-1}).
Written in the second order equation, it is expressed as $q'' = 6q^2 + z$.
Its Laurent series solution is given by
\begin{equation}
q(z) = (z-\alpha _0)^{-2} -\frac{\alpha _0}{10} (z-\alpha _0)^2 - \frac{1}{6}(z-\alpha _0)^3
 + \alpha _1 (z-\alpha _0)^4 +O((z-\alpha _0)^5),
\end{equation}
where $\alpha _1$ is an arbitrary parameter called a free parameter.
Counting from the lowest order $-2$, $\alpha _1$ is contained at order $6$
and this $6$ is called the Kovalevskaya exponent.
In general, for a Laurent series solution $x_i(z)$ of a quasi-homogeneous system,
if an arbitrary parameter is contained in the coefficient of $(z-\alpha _0)^{-a_i+j}$,
$j$ is called the Kovalevskaya exponent.

As an another example, let us consider the autonomous limit of the 
fourth order first Painlev\'{e} equation.
It is a four dimensional system defined by the following Hamiltonian
\begin{equation}
H_1(q_1, p_1, q_2, p_2) = 2p_1p_2 + 3p_2^2 q_1+q_1^4 - q_1^2q_2-q_2^2.
\end{equation}
The weight is given by $\deg (q_1,p_1,q_2,p_2) = (2,5,4,3)$.
The Hamiltonian vector field has two types of families of Laurent series solutions.
The one is starting from $q_1(z) = (z-\alpha _0)^{-2} + \cdots $,
and it contains three arbitrary parameters in the coefficients at order 2, 5 and 8.
Thus the Kovalevskaya exponents are $2,5,8$.
The other family starting from $q_1(z) = 3(z-\alpha _0)^{-2} + \cdots $
contains only two arbitrary parameters at order 8 and 10.
Thus the Kovalevskaya exponents are $8$ and $10$.

In general, for a given $m$-dimensional vector field $F$, if any solution is meromorphic,
the equation may have a family of solutions that contains $m$ arbitrary parameters determined by an initial condition.
In the above example $q_1(z) = (z-\alpha _0)^{-2} + \cdots $,
they are the pole $\alpha _0$ and three parameters in coefficients.
Such a family of Laurent series solutions is called the \textit{principal} Laurent series
that constructs an $m$ dimensional manifold $\mathcal{M}_m$.
A boundary of the manifold, if it exists, may be an $m-1$ dimensional manifold $\mathcal{M}_{m-1}$.
It is occupied by a different family of Laurent series solutions, that contains only $m-1$ arbitrary
parameters. Such a family is called \textit{non-principal} or \textit{lower}.
The latter family $q_1(z) = 3(z-\alpha _0)^{-2} + \cdots $ in the above example is this case.

Then, a natural question arises:
can we construct a lower Laurent series solutions from the principal one?
The purpose in this article is to consider this problem and give a systematic way
to obtain the lower one from the principal one.

For this purpose, we assume that for a given quasi-homogeneous vector field $F$,
there exists a quasi-homogeneous vector field $G$ that commutes with $F$.
Let us take an initial point on $\mathcal{M}_m$ and consider the solution $x(z)$
of $dx/dz = F$. Then, $x(z) \in \mathcal{M}_m$ for any $z\in \C$.
However, if we change the ``route" in the sense that 
$x(z)$ is governed by the flow of $G$ from some point,
a solution may reach at $\mathcal{M}_{m-1}$ along the orbit of $G$ 
and a principal Laurent series solutions degenerates to a lower one.
Indeed, for the above example, there is another Hamiltonian $H_2$ that commutes
with $H_1$ in Poisson bracket, see Example 4.13 for $H_2$. 
Based on this idea, we will obtain a lower Laurent series solutions from the principal one.

To be more precise, let $x(z_1; \alpha _0, \cdots ,\alpha _{m-1})$ be an $m$ parameter family of 
Laurent series solutions of the equation $dx/dz_1 = F$, 
where $\alpha _0, \cdots , \alpha _{m-1}$ are arbitrary parameters.
As was explained in the above example, the orders at which the parameters
appear in the Laurent expansion are called the Kovalevskaya exponents denoted by
$\kappa _0, \cdots ,\kappa _{m-1}$.
Suppose that there is a commuting vector field $G$ with $F$ that associates
the second equation $dx/dz_2 = G(x)$.
Assume that $\alpha _i=\alpha _i (z_2)$ depends on the second time $z_2$ so that 
$x(z_1, \alpha_0 (z_2), \cdots ,\alpha _{m-1}(z_2))$ solves both equations.
This assumption yields the differential equation of $(\alpha _0(z_2), \cdots ,\alpha _{m-1}(z_2))$
that is called the parameter flow.
The main results in this article state that under certain assumptions, the parameter flow
has the following properties (Prop.4.1, Thm.4.6 and Thm.4.12); 
\\[0.2cm]
\textbf{(i)} it is also quasi-homogeneous with the weight $(\kappa_0,\cdots ,\kappa_{m-1})$,
that is the Kovalevskaya exponents of the \textit{principal} Laurent series of $F$.
\\
\textbf{(ii)} the Kovalevskaya exponents of the Laurent series solution of the parameter
flow coincide with that of the \textit{lower} Laurent series of $F$.
\vskip.5\baselineskip

I this way, we can deduce the existence of the lower Laurent series solutions and its 
Kovalevskaya exponents of $F$ from the parameter flow.

For the list of weights, types of Laurent series solutions and their Kovalevskaya exponents
of four dimensional polynomial Painlev\'{e} equations,
the reader can refer to \cite{Chi4}.
The method developed in this article is applicable to all of them.


\section{Settings and the Kovalevskaya exponents}

Let $F=(f_1, \cdots  ,f_m)$ and $G = (g_1, \cdots, g_m)$ be
quasi-homogeneous polynomial vector fields on $\C^m$.
We consider the following partial differential equations
\begin{equation}
\frac{\partial x_i}{\partial z_1} = f_i(x), \quad 
\frac{\partial x_i}{\partial z_2} = g_i(x), \quad i=1,\cdots ,m,
\label{2-1}
\end{equation}
where $x = (x_1, \cdots  ,x_m)\in \C^m$ and $z_1, z_2 \in \C$.
We suppose the following.
\vskip.5\baselineskip

\textbf{(A1)} $F$ and $G$ are quasi-homogeneous: 
there exists a tuple of positive integers $(a_1, \cdots ,a_m) \in \N^m$
and $\gamma \in \N$ such that 
\begin{equation}
\left\{ \begin{array}{l}
\displaystyle f_i(\lambda ^{a_1}x_1 ,\cdots , \lambda ^{a_m}x_m) = \lambda ^{a_i+1}f_i(x_1, \cdots ,x_m)  \\[0.2cm]
\displaystyle g_i(\lambda ^{a_1}x_1 ,\cdots , \lambda ^{a_m}x_m) = \lambda ^{a_i+\gamma }g_i(x_1, \cdots ,x_m),  \\
\end{array} \right.
\label{A1}
\end{equation}
for any $\lambda \in \C$.
We call $\gamma $ the degree of $G$ with respect to the weight $(a_1, \cdots ,a_m)$.
The degree of $F$ is assumed to be $1$.
\vskip.5\baselineskip

\textbf{(A2)} $F$ and $G$ commute with each other with respect to the Lie bracket:
$[F, G] = 0$. This is equivalent to
\begin{equation}
\sum^m_{j=1} \left( f_j(x) \frac{\partial g_i}{\partial x_j}(x) 
  - g_j(x) \frac{\partial f_i}{\partial x_j}(x) \right) = 0, \quad i=1,\cdots ,m.
\label{A2}
\end{equation}

\textbf{(A3)} $F(x) = 0$ only when $x = 0$.
\vskip.5\baselineskip

In this section, we consider only the flow of $F$ and $z_1$ is denoted by $z$ for simplicity.
Let us consider the \textit{formal} series solution of $dx_i/dz = f_i(x)$ of the form
\begin{equation}
x_i(z) = c_i(z - \alpha _0)^{-q_i} + \sum^\infty_{j=1} d_{i,j} (z - \alpha _0)^{-q_i+j},
\label{2-4}
\end{equation}
where $q_i \in \N$, $\alpha _0$ is a possible singularity and $c_i, d_{i,j} \in \C$ are coefficients.
\vskip.5\baselineskip

\textbf{Theorem \thedef (\cite{Chi1}, Thm.2.9).}
Under the assumptions (A1) and (A3), any \textit{formal} Laurent series solution (\ref{2-4}) 
is a \textit{convergent} Laurent series solution of the form
\begin{equation}
x_i(z) = c_i(z - \alpha _0)^{-a_i} + \sum^\infty_{j=1} d_{i,j} (z - \alpha _0)^{-a_i+j}, 
\quad (c_1, \cdots ,c_m) \neq (0,\cdots ,0),
\label{2-5}
\end{equation}
where the exponents $(a_1, \cdots ,a_m)$ are the weight of $F$.
\vskip.5\baselineskip

This theorem means that there are no Laurent series solutions whose orders of poles are larger than $(a_1, \cdots , a_m)$.
Further, we can show that if $(c_1, \cdots , c_m) = (0,\cdots ,0)$, then $d_{i,j} = 0$ for $j=1, \cdots , a_{i}-1$.
This means that (\ref{2-5}) is a local holomorphic solution.
To prove it, the assumption (A3) is essentially used, while for the convergence of the series, (A1) is enough.
By substituting (\ref{2-5}) into the equation $dx_i/dz = f_i(x)$ and comparing
the coefficients of $(z-\alpha _0)^{j}$ in both sides, it turns out that 
$c_i$ and $d_{i,j}$ are given as follows.
\vskip.5\baselineskip

\textbf{Definition \thedef.}
A root $c = (c_1, \cdots ,c_m)\in \C^m,\,\, c\neq (0, \cdots ,0)$ of the equation
\begin{equation}
-a_ic_i = f_i(c_1, \cdots ,c_m), \quad i=1,\cdots ,m
\label{2-6}
\end{equation}
is called the indicial locus.
For a fixed indicial locus $c$, $d_{j} = (d_{1,j}, \cdots ,d_{m,j})$ is iteratively determined as a solution 
of the equation
\footnote{In this article, a column vector is often written as a row vector.}
\begin{equation}
(K(c) - j \cdot I)d_j = \text{( polynomial of $c$ and $d_k$ for $k=1, \cdots ,j-1$)},
\label{2-7}
\end{equation}
where $I$ is the identity matrix and $K(c)$ is defined by
\begin{equation}
K(c) = \left\{ \frac{\partial f_i}{\partial x_j}(c_1, \cdots ,c_m)
 + a_i \delta _{i,j} \right\}_{i,j=1}^m,
\end{equation}
that is called the Kovalevskaya matrix (K-matrix).
Its eigenvalues $\kappa_0, \kappa_1, \cdots ,\kappa_{m-1}$ are called the 
Kovalevskaya exponents (K-exponents) associated with $c$.
\vskip.5\baselineskip

The following results are well known, see \cite{Adl, Gor2}. 
\vskip.5\baselineskip

\textbf{Proposition \thedef.}
\\
\noindent (i)  $-1$ is always a Kovalevskaya exponent.
One of the eigenvectors of $K$ is given by $(a_1c_1, \cdots ,a_mc_m)$.
In what follows, we set $\kappa_0 = -1$.
\vskip.5\baselineskip

\noindent (ii) $\kappa =0$ is a Kovalevskaya exponent associated with $c$ if and only if $c$
is not an isolated root of the equation (\ref{2-6}).
\vskip.5\baselineskip

From (\ref{2-7}), it follows that 
if a positive integer $j$ is not an eigenvalue of $K(c)$, $d_j$ is uniquely determined.
If a positive integer $j$ is an eigenvalue of $K(c)$ and (\ref{2-7}) has a solution $d_j$, 
then $d_j + v$ is also a solution for any eigenvectors $v$.
This implies that the Laurent series solution (\ref{2-5}) contains an arbitrary parameter
(called a free parameter) in $d_j = (d_{1,j} ,\cdots ,d_{m,j})$.
Therefore, if (\ref{2-5}) represents a $k$-parameter family of Laurent series solutions
which contains $k-1$ free parameters other than $\alpha _0$, 
at least $k-1$ Kovalevskaya exponents have to be nonnegative integers and we need $k-1$ independent eigenvectors
associated with them.
Hence, the classical Painlev\'{e} test \cite{Gor2} for the necessary condition 
for the Painlev\'{e} property is stated as follows;
\vskip.5\baselineskip

\textbf{Classical Painlev\'{e} test.} 
If the system (\ref{2-1}) satisfying (A1) and (A3) has a $m$-parameter family 
of Laurent series solutions of the form (\ref{2-5}), there exists an indicial locus $c = (c_1,\cdots ,c_m)$ such that
all Kovalevskaya exponents except for $\kappa_0 =-1$ are nonnegative integers,
and the Kovalevskaya matrix is semisimple.
\vskip.5\baselineskip

In Chiba~\cite{Chi1}, the necessary and sufficient condition for the system (\ref{2-1}) to have 
an $m$-parameter family of Laurent series solutions is given, that is called the extended Painlev\'{e} test.
The gap of them is that even if the necessary condition for the classical Painlev\'{e} test is satisfied,
(\ref{2-7}) may not have a solution $d_j$.
In this case, the series solution (\ref{2-5}) is modified as a combination of powers of $(z-\alpha _0)$
and $\log (z-\alpha _0)$.
\vskip.5\baselineskip

\textbf{Definition \thedef.}
An indicial locus $c=(c_1, \cdots ,c_m)$ is called \textit{principal}
if the associated Laurent series solution (\ref{2-5}) exists and contains $m$ free parameters.
If the number of free parameters is smaller than $m$, the locus is 
called a \textit{lower} indicial locus.
\vskip.5\baselineskip

In the rest of this article, we assume that there exists an \textit{isolated} principal indicial locus $c$
of the vector field $F$; that is, 
all its K-exponents are positive integers other than $\kappa_0 = -1$
and the series solution (\ref{2-5}) contains $m$ free parameters (one of which is $\alpha _0$).
In this case, for each $\kappa_j\,\, (j \neq 0)$, we can take $d_{i,\kappa_j}$ as 
a free parameter for some $i$.
We denote it as $d_{i, \kappa_j} = \alpha _j$.
Then, all coefficients $d_{i,j}$ of (\ref{2-5}) are polynomials of $\alpha _1, \cdots ,\alpha _{m-1}$ and 
the solution is expressed as
\begin{align*}
x_i(z) &= c_i(z - \alpha _0)^{-a_i} + \sum^\infty_{j=1} d_{i,j} (z - \alpha _0)^{-a_i+j} \\
&= x_i(z; \alpha _0, \alpha _1, \cdots ,\alpha _{m-1}).
\end{align*}
The initial value of $x_i(z; \alpha _0, \alpha _1, \cdots ,\alpha _{m-1})$ is denoted by
\begin{align*}
& x_i(0; \alpha _0, \alpha _1, \cdots ,\alpha _{m-1})
= \Phi_i(\alpha _0, \alpha _1, \cdots ,\alpha _{m-1}) = \Phi_i(A),\\
& \Phi (A) = (\Phi_1(A), \cdots ,\Phi_m(A)),\quad A = (\alpha _0, \alpha _1, \cdots ,\alpha _{m-1}),
\end{align*}
which is well-defined for  small $|\alpha _0| \neq 0$ and $\Phi(A)$ is a locally biholomorphic map into $\C^m$.
In what follows, $d_{i,j}$ is denoted by $d_{i,j}(A)$ as a polynomial of $\alpha _0,\cdots ,\alpha _{m-1}$,
though it does not depend on $\alpha _0$ by the construction.
\vskip.5\baselineskip

\textbf{Proposition \thedef.} Put 
\[
\lambda \cdot A 
:= (\lambda ^{-1}\alpha _0, \lambda ^{\kappa_1}\alpha _1,\cdots ,\lambda ^{\kappa_{m-1}}\alpha _{m-1}).
\]
Functions $d_{i,j}$ and $\Phi_i$ are quasi-homogeneous satisfying
\begin{equation}
d_{i,j}(\lambda \cdot A) = \lambda ^j d_{i,j}(A), \quad 
\Phi_i(\lambda \cdot A) = \lambda ^{a_i} \Phi_i (A), \quad \lambda \in \C.
\label{2-9}
\end{equation}

\textbf{Proof.} Put $\widetilde{z} = \lambda ^{-1}z$ and $\widetilde{x}_i = \lambda ^{a_i}x_i$.
Then, $\widetilde{x}(\widetilde{z})$ satisfies the same equation as $x(z)$ because of (\ref{A1}).
Let 
\begin{equation}
x_i = c_i(z-\alpha _0)^{-a_i} + \sum^\infty_{j=1}d_{i,j}(A) (z-\alpha _0)^{-a_i+j}
\label{2-10}
\end{equation}
be a Laurent series solution with free parameters $\alpha _0, \cdots ,\alpha _{m-1}$.
Then
\begin{align}
& \Phi_i(A) = c_i(-\alpha _0)^{-a_i} + \sum^\infty_{j=1}d_{i,j}(A) (-\alpha _0)^{-a_i+j}, 
\label{2-11} \\
& \Phi_i(\lambda \cdot A) 
= \lambda ^{a_i}\left(  c_i(-\alpha _0)^{-a_i} 
   + \sum^\infty_{j=1}\lambda ^{-j} d_{i,j}(\lambda \cdot A) (-\alpha _0)^{-a_i+j} \right).
\label{2-12}
\end{align}
Similarly, consider the Laurent series solution of $\widetilde{x}$, whose locus $c$ is the same as that of (\ref{2-11}):
\begin{equation}
\widetilde{x}_i = c_i(\widetilde{z}-\widetilde{\alpha} _0)^{-a_i} 
+ \sum^\infty_{j=1}d_{i,j}(\widetilde{A}) (\widetilde{z}-\widetilde{\alpha} _0)^{-a_i+j},
\quad \widetilde{A}=(\widetilde{\alpha }_0, \cdots , \widetilde{\alpha }_{m-1}).
\label{2-13}
\end{equation}
Since $x(z)$ and $\widetilde{x}(\widetilde{z})$ satisfy the same equation, $d_{i,j}$ in (\ref{2-10})
and (\ref{2-13}) are the common function of $A$, though we can choose different values of free parameters.
If we put $\widetilde{\alpha }_0 = \lambda ^{-1}\alpha _0$, (\ref{2-13}) is rewritten as
\begin{align}
& \lambda ^{a_i} x_i = c_i\lambda ^{a_i} (z-\alpha _0)^{-a_i}
 + \sum^\infty_{j=1}d_{i,j}(\widetilde{A}) \lambda ^{a_i-j} (z-\alpha _0)^{-a_i+j}, \nonumber \\
\Rightarrow \quad  & \Phi_i(A) = c_i (-\alpha _0)^{-a_i}
 + \sum^\infty_{j=1}d_{i,j}(\widetilde{A}) \lambda ^{-j} (-\alpha _0)^{-a_i+j}.
\label{2-14}
\end{align}
It follows from (\ref{2-11}) and (\ref{2-14}) that $d_{i,j}(A) = \lambda ^{-j}d_{i,j}(\widetilde{A})$.
When $\kappa_j$ is one of the K-exponents, $d_{i,\kappa_j}(A) = \alpha _j$ for some $i$ by the definition.
Hence, we have $\alpha _j = \lambda ^{-\kappa_j} \widetilde{\alpha }_j $.
This shows $\widetilde{A} = \lambda \cdot A$ and $d_{i,j}(A) = \lambda ^{-j}d_{i,j}(\lambda \cdot A)$.
Therefore, (\ref{2-11}) and (\ref{2-12}) prove the desired result. $\hfill \Box$
\vskip.5\baselineskip

\textbf{Corollary \thedef.}
$d_{i,j}(A) \neq 0$ only when there exists a tuple of integers $(n_1,\cdots ,n_{m-1}) \neq (0,\cdots ,0)$
such that $n_1 \kappa_1 + \cdots + n_{m-1}\kappa_{m-1} = j$.
\vskip.5\baselineskip
\textbf{Proof.} 
Assume that a monomial $\alpha _1^{n_1}\alpha _2^{n_2} \cdot \cdots \cdot \alpha _{m-1}^{n_{m-1}}$
is contained in $d_{i,j}(A)$.
Substituting it into (\ref{2-9}) proves the desired result. $\hfill \Box$
\vskip.5\baselineskip

For a quasi-homogeneous vector field $dx/dz = F(x)$ satisfying (\ref{A1}), 
let $x = \varphi (y_1, \cdots ,y_m),\,\, \varphi = (\varphi _1, \cdots ,\varphi _m)$ be 
a (locally) holomorphic coordinate transformation satisfying
\begin{equation*}
\varphi _i (\lambda ^{q_1}y_1, \cdots ,\lambda ^{q_m} y_m) = \lambda ^{a_i}\varphi_i (y_1, \cdots ,y_m),
\quad i=1,\cdots ,m,
\end{equation*}
where $(q_1, \cdots ,q_m) \in \Z^m$ is an arbitrary tuple of integers and $(a_1, \cdots ,a_m)$ is the same as in (\ref{A1}).
Remark that (\ref{2-9}) is just in this case.
By the transformation, $dx/dz = F(x)$ is transformed into the system
\begin{align*}
\frac{dy}{dz} = (D\varphi )^{-1}F(\varphi (y)) := \widetilde{F}(y),
\end{align*}
where $D\varphi $ is the Jacobi matrix.
\vskip.5\baselineskip

\textbf{Theorem \thedef \, (\cite{Chi1}, Thm.2.5).}

$\widetilde{F}(y)$ is quasi-homogeneous with respect to the weight $(q_1, \cdots ,q_m)$ whose degree is the same 
as that of $F$.
If $c$ is an indicial locus of $F$, $\widetilde{c}= \varphi ^{-1}(c)$ is an indicial locus of $\widetilde{F}$.
The K-exponents of $\widetilde{F}$ at $\widetilde{c}$ are the same as those of $F$ at $c$.
\vskip.5\baselineskip

\textbf{Example \thedef.} Consider the 2-dim system
\begin{align*}
\frac{dx}{dz} = y, \quad \frac{dy}{dz} = 6 x^2.
\end{align*}
This satisfies the assumptions (A1) and (A3) with the weight $(a_1, a_2) = (2,3)$.
The indicial locus is uniquely given by $(c_1, c_2) = (1,-2)$.
Thus, the Laurent series solution starts from $(x, y) = (T^{-2},\, -2 T^{-3}),\,\, T = z-\alpha _0$.
The K-exponent associated with the locus is $\kappa_1 = 6$.
Hence, a free parameter appears in $d_{1,6}$ and/or $d_{2,6}$.
Indeed, we can verify that $d_{2,6} = 4d_{1,6}$ and a solution is given by
\begin{align}
\left(
\begin{array}{@{\,}c@{\,}}
x \\
y
\end{array}
\right)&=\left(
\begin{array}{@{\,}c@{\,}}
0 \\
-2
\end{array}
\right) T^{-3} +\left(
\begin{array}{@{\,}c@{\,}}
1 \\
0
\end{array}
\right) T^{-2}+\left(
\begin{array}{@{\,}c@{\,}}
0 \\
4 d_{1,6}
\end{array}
\right) T^{3}+\left(
\begin{array}{@{\,}c@{\,}}
d_{1,6} \\
0
\end{array}
\right) T^{4} \nonumber \\
&\qquad + \left(
\begin{array}{@{\,}c@{\,}}
0 \\
10 d_{1,6}^2/13
\end{array}
\right) T^{9} + \left(
\begin{array}{@{\,}c@{\,}}
d_{1,6}^2/13 \\
0
\end{array}
\right) T^{10} + \cdots .
\label{ex1}
\end{align}
In this case, we put $\alpha _1 : = d_{1,6}$.
All other coefficients are polynomials of $\alpha _1$
and we can confirm Prop. 2.5.

\section{Properties of the vector field $G$}

To consider the vector field $G = (g_1, \cdots ,g_m)$, we prepare several formulae and notations.
The derivative of (\ref{A1}) at $\lambda =1$ yields
\begin{align}
\sum^m_{j=1}a_jx_j \frac{\partial f_i}{\partial x_j}(x) = (a_i+1)f_i(x), \quad
\sum^m_{j=1}a_jx_j \frac{\partial g_i}{\partial x_j}(x) = (a_i+\gamma )g_i(x).
\label{3-1}
\end{align}
Putting  $x=c$ to the latter one, we obtain
\begin{align}
-\sum^m_{j=1} f_j(c) \frac{\partial g_i}{\partial x_j}(c) = (a_i+\gamma )g_i(c).
\label{3-2}
\end{align}
The derivative of (\ref{A1}) by $x_j$ gives
\begin{equation}
\left\{ \begin{array}{l}
\displaystyle \frac{\partial f_i}{\partial x_j}(\lambda ^{a_1}x_1,\cdots ,\lambda ^{a_m}x_m)
 = \lambda ^{a_i+1-a_j}\frac{\partial f_i}{\partial x_j}(x_1, \cdots ,x_m),  \\[0.3cm]
\displaystyle \frac{\partial g_i}{\partial x_j}(\lambda ^{a_1}x_1,\cdots ,\lambda ^{a_m}x_m)
 = \lambda ^{a_i+\gamma -a_j}\frac{\partial g_i}{\partial x_j}(x_1, \cdots ,x_m).
\end{array} \right.
\label{3-3}
\end{equation}
As before, we assume that $c$ is a principal locus so that there is a Laurent series solution
\begin{equation}
x_i(z) = (z -\alpha _0)^{-a_i} \left( c_i + \sum^\infty_{k=1} d_{i,k}(A) (z-\alpha _0)^k \right)
 =: (z -\alpha _0)^{-a_i}y_i(z)
\label{3-4}
\end{equation}
including $m$ free parameters $A = (\alpha _0, \cdots ,\alpha _{m-1})$.
Substituting this solution into $g_i(x)$ gives
\begin{align}
g_i(x(z)) &= (z-\alpha _0)^{-a_i-\gamma } g_i(y_1, \cdots ,y_m) \nonumber \\
& =  (z-\alpha _0)^{-a_i-\gamma } \sum^\infty_{k=0}g_{i,k}(A) (z-\alpha _0)^k,
\label{3-5}
\end{align}
where $g_{i,k}(A)$ is the coefficient of the Taylor expansion of $g_i(y_1, \cdots ,y_m)$ in $z-\alpha _0$.
We denote $G_k(A) := (g_{1,k}(A), \cdots ,g_{m,k}(A))$.
In particular $g_{i,0} = g_i(c)$ and $G_0 = G(c)$.
\vskip.5\baselineskip

\textbf{Proposition \thedef.} 
For a given indicial locus $c$, the identity $(K(c) + \gamma ) G(c) = 0$ holds.
In particular, if $-\gamma $ is not a K-exponent, then $G(c) = 0$
(here, we need not assume that $c$ is principal).
\vskip.5\baselineskip

\textbf{Proof.}
(\ref{A2}) and (\ref{3-2}) provide
\begin{align*}
0 = \sum^m_{j=1}\left( f_j(c) \frac{\partial g_i}{\partial x_j}(c) 
     - g_j(c) \frac{\partial f_i}{\partial x_j}(c)  \right) 
= - \sum^m_{j=1} \frac{\partial f_i}{\partial x_j}(c) g_j(c) - (a_i + \gamma ) g_i(c),
\end{align*}
which proves the proposition. $\hfill \Box$
\vskip.5\baselineskip

\textbf{Corollary \thedef.}
If $G(c) \neq 0$, then $-\gamma $ is a K-exponent.
In particular, when $\gamma \neq 1$, there exists a lower indicial locus.
\vskip.5\baselineskip

\textbf{Corollary \thedef.} Suppose $c$ is a principal indicial locus.
\\
(i) If $\gamma \geq 2$, $G(c) = 0$.
\\
(ii) If $\gamma =1$, $G(c) = k (a_1c_1, \cdots ,a_m c_m)$ for some $k\in \C$.
\vskip.5\baselineskip

\textbf{Proof.} (i) By the assumption, there are no negative K-exponents other than $\kappa_0 = -1$.
(ii) When $\gamma =1$, $G(c) = 0$ or $G(c)$ is an eigenvector of $\kappa_0 = -1$.
 Since $\kappa_0 = -1$ is a simple eigenvalue by the assumption, the statement (ii) follows from Prop.2.3. 
$\hfill \Box$
\vskip.5\baselineskip

More generally, the next theorem holds.
\vskip.5\baselineskip

\textbf{Theorem \thedef.} Let $c$ be a principal indicial locus.
\vskip.5\baselineskip

The equality $(K(c)+\gamma -k)G_k(A) = 0$ holds for $k=0,1,\cdots ,\gamma -1$.
\vskip.5\baselineskip

\textbf{Corollary \thedef.} When $\gamma \geq 2$,
$G_0 = G_1= \cdots  = G_{\gamma -2}= 0$ and $G_{\gamma -1}$ is of the form
$G_{\gamma -1} = h (A) \cdot (a_1c_1, \cdots ,a_mc_m)$, where $h(A)$ is a certain polynomial
of $\alpha _1 , \cdots ,\alpha _{m-1}$.
\vskip.5\baselineskip

\textbf{Proof.}
The case $\gamma =1$ (i.e. $k=0$ in the theorem) had been proved in Prop.3.1.
Thus we consider $\gamma \geq 2$.

Since $x(z)$ in (\ref{3-4}) is a solution of $dx/dz = f(x)$,
\begin{align*}
f_j(x(z)) &= -a_jc_j(z-\alpha _0)^{-a_j-1} - \sum^\infty_{k=1} (a_j-k)d_{j,k}(z-\alpha _0)^{-a_j+k-1} \\
&= -(z-\alpha _0)^{-1}a_jx_j(z) + \sum^\infty_{k=1} k d_{j,k} (z-\alpha _0)^{-a_j+k-1}.
\end{align*}
Eq.(\ref{3-3}) shows
\begin{align*}
\frac{\partial f_i}{\partial x_j}(x(z)) = (z-\alpha _0)^{a_j-a_i-1}\frac{\partial f_i}{\partial x_j}(y),
\quad \frac{\partial g_i}{\partial x_j}(x(z)) = (z-\alpha _0)^{a_j-a_i-\gamma }\frac{\partial g_i}{\partial x_j}(y),
\end{align*}
where $y = (y_1, \cdots ,y_m)$ is defined in (\ref{3-4}).
Substituting them into (\ref{A2}) with (\ref{3-1}) and (\ref{3-5}), we have
\begin{align*}
0 &= \sum^m_{j=1}\left( (z-\alpha _0)^{a_j-a_i-1} \frac{\partial f_i}{\partial x_j}(y) \cdot
(z-\alpha _0)^{-a_j-\gamma }\sum^\infty_{k=0} g_{j,k}(A)(z-\alpha _0)^{k} \right) \\
& \qquad + \sum^m_{j=1}  \frac{\partial g_i}{\partial x_j}(x(z))  
\left(  a_jx_j(z) (z-\alpha _0)^{-1} 
   - \sum^\infty_{k=1}k d_{j,k} (z-\alpha _0)^{-a_j+k-1} \right) \\
&= \sum^m_{j=1} (z-\alpha _0)^{-a_i-\gamma -1}\frac{\partial f_i}{\partial x_j}(y) 
   \sum^\infty_{k=0}g_{j,k}(A) (z-\alpha _0)^{k} \\
& \quad + (a_i + \gamma )  (z-\alpha _0)^{-a_i-\gamma -1}\sum^\infty_{k=0} g_{i,k}(A) (z-\alpha _0)^{k} \\
& \qquad -\sum^m_{j=1} (z-\alpha _0)^{a_j-a_i-\gamma } 
    \frac{\partial g_i}{\partial x_j}(y)\sum^\infty_{k=1}kd_{j,k} (z-\alpha _0)^{-a_j+k-1}.
\end{align*}
Multiplied by $(z-\alpha _0)^{a_i+\gamma +1}$, this is rewritten as  
\begin{equation}
\sum^\infty_{k=0} (z-\alpha _0)^{k}\sum^m_{j=1} 
   \left(  \frac{\partial f_i}{\partial x_j}(y) + (a_j + \gamma ) \delta _{i,j} \right) g_{j,k}(A)
-\sum^m_{j=1} \frac{\partial g_i}{\partial x_j}(y)\sum^\infty_{k=1}kd_{j,k} (z-\alpha _0)^{k}=0.
\label{3-6}
\end{equation}
To estimate the last term, we use induction. Recall that we consider $\gamma \geq 2$.
Assume that $(K(c) +\gamma -k)G_k(A) = 0$ holds for $k=0,1,\cdots ,n-1$, where $n\leq \gamma -1$.
Then $G_0 = G_1 = \cdots = G_{n-1}= 0$ holds because $-\gamma , -\gamma +1, \cdots , -\gamma + n -1\, (\leq -2)$
are not K-exponents.
Hence, (\ref{3-6}) gives 
\begin{equation*}
\sum^\infty_{k=n} (z-\alpha _0)^{k}\sum^m_{j=1} 
    \left(  \frac{\partial f_i}{\partial x_j}(y) + (a_j + \gamma ) \delta _{i,j} \right) g_{j,k}(A)
-\sum^m_{j=1} \frac{\partial g_i}{\partial x_j}(y)\sum^\infty_{k=1}kd_{j,k} (z-\alpha _0)^{k}=0
\end{equation*}
(the first summation starts from $k=n$). Divide  by $(z-\alpha _0)^n$ and consider the limit $z\to \alpha _0$:
\begin{equation}
\sum^m_{j=1} \left(  \frac{\partial f_i}{\partial x_j}(c) + (a_j + \gamma ) \delta _{i,j} \right) g_{j,n}(A)
-\lim_{z\to \alpha _0} \sum^m_{j=1} \frac{\partial g_i}{\partial x_j}(y)\sum^\infty_{k=1}kd_{j,k}
\frac{(z-\alpha _0)^k}{(z-\alpha _0)^n}=0.
\label{3-7}
\end{equation}
By the definition of $g_{i,k}(A)$, we have
\begin{align*}
 & g_i(y_1, \cdots ,y_m) = \sum^\infty_{k=0} g_{i,k}(A)(z-\alpha _0)^k =\sum^\infty_{k=n} g_{i,k}(A)(z-\alpha _0)^k, \\
 & y_j = c_j + \sum^\infty_{k=1} d_{j,k} (z-\alpha _0)^k.
\end{align*}
The derivation of both sides by $z$ yields
\begin{align*}
& \sum^m_{j=1} \frac{\partial g_i}{\partial x_j}(y_1,\cdots ,y_m)\sum^\infty_{k=1}k d_{j,k} (z-\alpha _0)^{k-1}
 = \sum^\infty_{k=n} k g_{i,k}(A) (z-\alpha _0)^{k-1} \\
\Rightarrow & \sum^m_{j=1} \frac{\partial g_i}{\partial x_j}(y_1,\cdots ,y_m)\sum^\infty_{k=1}k d_{j,k} 
\frac{(z-\alpha _0)^{k-1}}{(z-\alpha _0)^{n-1}}
 = ng_{i,n}(A) + O(z-\alpha _0). 
\end{align*}
As $z\to \alpha _0$,
\begin{equation*}
\lim_{z\to \alpha _0}\sum^m_{j=1} \frac{\partial g_i}{\partial x_j}(y_1,\cdots ,y_m)\sum^\infty_{k=1}k d_{j,k} 
\frac{(z-\alpha _0)^{k}}{(z-\alpha _0)^{n}}
 = ng_{i,n}(A)
\end{equation*}
This and (\ref{3-7}) gives
\begin{equation}
\sum^m_{j=1} \left(  \frac{\partial f_i}{\partial x_j}(c) + (a_j + \gamma -n) \delta _{i,j} \right) g_{j,n}(A)=0.
\end{equation}
This proves that $(K(c) +\gamma -k)G_k(A) = 0$ holds for $k=n$.
The induction step continues up to $k=\gamma -1$. $\hfill \Box$
\vskip.5\baselineskip

Even when the degree $\gamma $ of a given quasi-homogeneous equation is larger than $1$,
the K-exponents are defined in a similar manner as follows.
Let us consider $dx_i/dz_2 = g_i(x)$ given in (\ref{A1}) having the degree $\gamma $.
We consider the Puiseux series solution of the form
\begin{equation}
x_i(z_2) = p_i (z_2-\beta_0)^{-a_i/\gamma }
 + \sum^\infty_{k=1} q_{i,k}(z_2-\beta_0)^{(-a_i+k)/\gamma }, \quad i=1,\cdots ,m,
\end{equation}
where $\beta_0$ is a singularity and $p_i,\, q_{i,k}$ are constants to be determined.
Substituting it into the equation, it turns out that an indicial locus $p=(p_1, \cdots ,p_m)$
 is given as a root the equation
\begin{equation}
-\frac{a_i}{\gamma }p_i = g_i(p_1, \cdots ,p_m), \quad i=1,\cdots ,m.
\label{3-10}
\end{equation}
For a given $p$, $q_j = (q_{1,j}, \cdots ,q_{m,j})$ is iteratively determined as a solution of
\begin{equation}
(K_\gamma (p) - \frac{j}{\gamma }\cdot I ) q_j = \text{(polynomial of $p$ and $q_{k}$ for $k=1,\cdots ,j-1$)},
\label{3-11}
\end{equation}
where the K-matrix $K_\gamma (p)$ is defined by
\begin{equation}
K_\gamma (p) = \left\{ \frac{\partial g_i}{\partial x_j}(p) + \frac{a_i}{\gamma } \delta _{i,j} \right\}_{i,j=1}^m.
\label{3-12}
\end{equation}
The eigenvalues $\rho_0, \rho_1, \cdots  ,\rho_{m-1}$ are called the K-exponents.
If $j/\gamma $ is a K-exponent, $q_j$ contains a free parameter.
As in Prop.2.3, $\rho_0 = -1$ is always a K-exponent with the eigenvector $(a_1p_1, \cdots , a_mp_m)$.
When $\rho_1, \cdots , \rho_{m-1} \in \N/\gamma $, $p$ is called a principal indicial locus.

\section{Flow of the free parameters}

Let $\varphi ^F_{z_1}$ and $\varphi ^G_{z_2}$ be the flow of $F$ and $G$, respectively;
$\varphi ^F_{z_1}$ maps $x(0)$ to $x(z_1)$ along the orbit of vector field $F$, and similarly for $\varphi ^G_{z_2}$.
To obtain a solution of the system (\ref{2-1}) as a function of $z_1$ and $z_2$,
for a solution $x(z_1; A) = \varphi ^F_{z_1} \circ \Phi (A)$ of $\partial x/\partial z_1 = F(x)$, 
we assume that $A = A(z_2)$ is a function of $z_2$.
Let us substitute this $x$ into the second equation $\partial x / \partial z_2 = G(x)$.
\begin{equation*}
\text{(left hand side)} = \frac{\partial \varphi ^F_{z_1}}{\partial x}(\Phi (A))\frac{d\Phi}{dz_2}(A)
= \frac{\partial \varphi ^F_{z_1}}{\partial x}(\Phi (A)) 
\sum^{m-1}_{l=0} \frac{\partial \Phi}{\partial \alpha _l}(A)\frac{d\alpha _l}{dz_2}.
\end{equation*}
It is known that (\ref{A2}) is equivalent to the identity 
$\varphi ^F_t \circ \varphi ^G_s = \varphi ^G_s \circ \varphi ^F_t$.
The derivative of it at $s = 0$ gives
\begin{equation*}
\frac{\partial \varphi ^F_t}{\partial x}(x) G(x) = G(\varphi ^F_t (x)).
\end{equation*}
Hence, 
\begin{equation*}
\text{(right hand side)} = G(x(z_1; A)) = G(\varphi ^F_{z_1} \circ \Phi (A))
 = \frac{\partial \varphi ^F_{z_1}}{\partial x} (\Phi (A)) G(\Phi (A)).
\end{equation*}
This proves
\begin{equation}
\sum^{m-1}_{l=0} \frac{\partial \Phi}{\partial \alpha _l}(A)\frac{d\alpha _l}{d z_2} = G(\Phi (A))
\quad \Rightarrow \quad \frac{dA}{dz_2} = \left( \frac{\partial \Phi}{ \partial A}(A)\right)^{-1}G(\Phi (A)),
\label{4-1}
\end{equation}
that gives the flow of $A(z_2)$.
\vskip.5\baselineskip

\textbf{Proposition \thedef.} The right hand side of (\ref{4-1}) is \\
(i) independent of $\alpha _0$ and polynomial in $\alpha _1, \cdots ,\alpha _{m-1}$, \\
(ii) quasi-homogeneous of degree $\gamma $ with respect to the weight $(\kappa_0, \cdots , \kappa_{m-1})$.
\vskip.5\baselineskip

\textbf{Proof.}
(i) Denote $c_i = d_{i,0}$ for simplicity to express (\ref{2-5}) as 
$x_i(z_1) = \sum^\infty_{j=0}d_{i,j}(z_1-\alpha _0)^{-a_i + j}$.
By substituting it to $\partial x_i/\partial z_2 = g_i(x)$, we obtain
\begin{align*}
\text{(left hand side)} &= \sum^\infty_{j=0} (a_i-j)d_{i,j}(z_1-\alpha _0)^{-a_i+j-1}\frac{d\alpha _0}{dz_2}
   + \sum^\infty_{j=0} \frac{dd_{i,j}}{dz_2}(z_1-\alpha _0)^{-a_i+j}, \\
\text{(right hand side)} &= 
(z_1-\alpha _0)^{-a_i-\gamma } g_i(\sum^\infty_{j=0} d_{1,j} (z_1-\alpha _0)^j, 
    \cdots , \sum^\infty_{j=0} d_{m,j} (z_1-\alpha _0)^j), \\
&= (z_1-\alpha _0)^{-a_i-\gamma } \sum^\infty_{j=0}g_{i,j}(A)(\alpha _1, \cdots ,\alpha _{m-1}) \cdot (z_1-\alpha _0)^j,
\end{align*}
where $g_{i,j}(A)$ is a $j$-th coefficient of the Taylor expansion of 
$g_i(\sum^\infty_{j=0} d_{1,j} (z_1-\alpha _0)^j, \cdots )$ in $z_1-\alpha _0$, that was introduced in Sec.3.
This is a polynomial of $\alpha _1, \cdots ,\alpha _{m-1}$ because so is $d_{i,j}$.
Comparing the coefficients of $(z_1-\alpha _0)^{-a_i-1}$ in both sides, we obtain
\begin{equation}
a_ic_i \frac{d\alpha _0}{dz_2} = g_{i, \gamma -1}(\alpha _1,\cdots ,\alpha _{m-1}).
\label{4-2}
\end{equation}
Similarly, coefficients of $(z_1-\alpha _0)^{-a_i + \kappa_l}$ provide
\begin{equation}
(a_i-1-\kappa_l)d_{i,1+\kappa_l} \frac{d\alpha _0}{dz_2}+\frac{d\alpha _l}{dz_2}
 = g_{i, \kappa_l + \gamma }(\alpha _1,\cdots ,\alpha _{m-1}),\quad l=1,\cdots ,m-1,
\label{4-3}
\end{equation}
where $i$ is chosen so that $\alpha _l = d_{i,\kappa_l}$.
Eq.(\ref{4-2}) shows that $d\alpha _0/dz_2$ is independent of $\alpha _0$.
Thus, (\ref{4-3}) proves that $d\alpha _l/dz_2$ is also independent of $\alpha _0$
and is polynomial in $\alpha _1, \cdots ,\alpha _{m-1}$ for $l=1,\cdots ,m-1$.

(ii) Eq.(\ref{4-2}) and (\ref{4-3}) provide another expression of (\ref{4-1}) as
\begin{equation}
\left\{ \begin{array}{l}
\displaystyle \frac{d\alpha _0}{dz_2} = \frac{g_{i,\gamma -1}(A)}{a_i c_i} =: \hat{g}_0(A),  \\[0.3cm]
\displaystyle \frac{d\alpha _l}{dz_2} = g_{i,\kappa_l+\gamma } (A)
 - (a_i-1-\kappa_l) d_{i,1+\kappa_l}(A) \frac{g_{i,\gamma -1}(A)}{a_i c_i} =:\hat{g}_l(A).
\end{array} \right.
\label{4-4}
\end{equation}
Remark that $g_{i,\gamma -1}(A)/a_i c_i$ is independent of $i$ and well-defined even when $c_i=0$
because of Corollary 3.2 and 3.4 ($\hat{g}_0(A)$ here is $h(A)$ there).
As above, we put
\begin{equation*}
g_i(\sum^\infty_{j=0} d_{1,j}(A) (z_1-\alpha _0)^j, 
    \cdots , \sum^\infty_{j=0} d_{m,j}(A) (z_1-\alpha _0)^j)
 = \sum^\infty_{j=0}g_{i,j}(A) \cdot (z_1-\alpha _0)^j.
\end{equation*}
Since $\lambda ^j d_{i,j}(A)= d_{i,j}(\lambda \cdot A)$ by Prop.2.5,
the left hand side above is invariant by $A \mapsto \lambda \cdot A,\, 
\alpha _0\mapsto \lambda ^{-1}\alpha _0$ and $z_1\mapsto \lambda ^{-1}z_1$
(recall that 
$\lambda \cdot A = (\lambda ^{-1}\alpha _0, \lambda^{\kappa_1} \alpha _1, 
\cdots , \lambda^{\kappa_{m-1}} \alpha _{m-1})$).
Thus, we obtain $g_{i,j}(\lambda \cdot A) = \lambda ^j g_{i,j}(A)$ from the right hand side.
This proves
\begin{equation}
\hat{g}_l (\lambda^{\kappa_1} \alpha _1, \cdots , \lambda^{\kappa_{m-1}} \alpha _{m-1})
   = \lambda ^{\kappa_l+\gamma }\hat{g}_l(\alpha _1, \cdots ,\alpha _{m-1})
\label{4-5}
\end{equation}
for $l=0, 1,\cdots ,m-1$. $\hfill \Box$
\vskip.5\baselineskip

Now we investigate the system (\ref{4-1}) or equivalently (\ref{4-4});
\begin{equation}
\frac{dA}{dz_2} = \left( \frac{\partial \Phi}{ \partial A}(A)\right)^{-1}G(\Phi (A))
\quad \Leftrightarrow \quad
\left\{ \begin{array}{ll}
\displaystyle \frac{d\alpha _0}{dz_2} = \hat{g}_0(A)  \\[0.2cm]
\displaystyle \frac{d\alpha _i}{dz_2} = \hat{g}_i(A), \quad i=1,\cdots ,m-1 \\
\end{array} \right.
\label{4-6}
\end{equation}
satisfying (\ref{4-5}).
The next lemma immediately follows from (\ref{4-5}).
\vskip.5\baselineskip

\textbf{Lemma \thedef.}
A monomial $\alpha _1^{n_1} \alpha _2^{n_2} \cdots \alpha _{m-1}^{n_{m-1}}$
can be contained in $\hat{g}_l(A)$ only when a tuple of nonnegative integers 
$(n_1, \cdots ,n_{m-1})$ satisfies 
$\kappa _1 n_1 + \cdots  + \kappa _{m-1}n_{m-1} = \kappa_l + \gamma$.
\vskip.5\baselineskip

From the lemma, it turns out that $\hat{g}_l(A)$ for $l=1,\cdots ,m-1$ does not contain
a constant term when $c$ is an isolated principal indicial locus (i.e. $\kappa_l > 0$).
Applying the lemma to $l=0$ gives $\kappa _1 n_1 + \cdots  + \kappa _{m-1}n_{m-1} = \gamma - 1$.
When $\gamma =1$, $\hat{g}_0$ is a constant function.
When $\gamma \geq 2$, $\hat{g}_0$ does not contain a constant term.
For example when $\gamma =2$, there exists a K-exponent $\kappa =1$.
More generally, there exists a K-exponent smaller than $\gamma $ if $\hat{g}_0$ is not identically zero.

Let $\xi = (\xi_0 ,\cdots ,\xi_{m-1})$ be an indicial locus of (\ref{4-6}) and 
$\rho_0 = -1,\, \rho_1 ,\cdots ,\rho_{m-1}$ its K-exponents.
\vskip.5\baselineskip

\textbf{Lemma \thedef.} 
Besides $\rho_0 = -1$, (\ref{4-6}) has another negative exponent $\rho_1 = -1/\gamma $.
\vskip.5\baselineskip

\textbf{Proof.}
The proof is based on the fact that the right hand side of (\ref{4-6}) is independent of $\alpha _0$.
The K-matrix of (\ref{4-6}) defined by (\ref{3-12}) is given as
\begin{equation}
K_\gamma (\xi) = \left(
\begin{array}{@{\,}c|ccc@{\,}}
0 &\partial \hat{g}_0/\partial \alpha _1 \!\!\!\!\!\!\!\!\!\!& \cdots &\!\!\!\!\!\!\!\!\!\! 
    \partial \hat{g}_0/\partial \alpha _{m-1} \\ \hline
 & & & \\
0 & & \quad \displaystyle \left\{ \frac{\partial \hat{g}_i}{\partial \alpha _j} \right\}_{i,j=1}^{m-1} &\\
 & & & 
\end{array}
\right) + \left(
\begin{array}{@{\,}c|ccc@{\,}}
-1/\gamma & & 0 & \\ \hline
& \kappa _1/\gamma   & & \\
0 & &\!\! \ddots & \\
& & &\!\! \kappa _{m-1}/\gamma 
\end{array}
\right).
\label{4-7}
\end{equation}
Hence, its eigenvalues are $-1/\gamma $ and K-exponents of the subsystem 
$d\alpha _i/dz_2 = \hat{g}_i (A),\,\, i=1,\cdots ,m-1$.
Since the subsystem has a K-exponent $-1$, the proof is completed. $\hfill \Box$
\vskip.5\baselineskip

\textbf{Remark \thedef.}
The first expression of (\ref{4-6}) implies that it is obtained from $dx/dz_2 = G(x)$
by the coordinate transformation $x = \Phi (A)$.
For the system $dx/dz_2 = G(x)$, let $p = (p_1, \cdots , p_m)$ be an indicial locus.
Then, $\xi =  \Phi^{-1}(p)$ is an indicial locus of (\ref{4-6}) as long as 
$\xi$ is in the domain on which $\Phi$ is a diffeomorphism.
In this case, the K-exponents of (\ref{4-6}) at $\Phi^{-1}(p)$ coincide with those of $G$ at $p$
due to Prop.2.5 and Thm.2.7.


\subsection{$\gamma =1$}

When $\gamma = 1$, $d\alpha _0/dz_2 = \hat{g}_0(A) = g_i(c)/(a_ic_i) \,\, (=:\text{constant}\, k_1)$ due to (\ref{4-2}).
It is solved as $\alpha _0 = k_1 z_2 + k_2$.
Thus, we consider the system
\begin{equation}
\frac{d\alpha _i}{dz_2} = \hat{g}_i (\alpha _1,\cdots ,\alpha _{m-1}), \quad i=1,\cdots ,m-1,
\label{4-8}
\end{equation}
which is of degree $\gamma=1$ with respect to the weight $(\kappa _1,\cdots , \kappa _{m-1})$.
Let $\xi = (\xi_1,\cdots ,\xi_{m-1})$ be an indicial locus 
and $\rho_0=-1, \rho_2, \cdots ,\rho_{m-1}$ its K-exponents 
(we skip $\rho_1=-1/\gamma $ shown in Lemma 4.3 because it arises from the equation of $\alpha _0$).
The Laurent series solution is written as
\begin{equation}
\alpha _i(z_2) = (z_2- \beta_0)^{-\kappa _i} 
   \left( \xi_i + \sum^\infty_{j=1} \eta_{i,j} (z_2-\beta_0)^{j}\right)
 =:  (z_2- \beta_0)^{-\kappa _i}y_i,
\label{4-9}
\end{equation}
with coefficients $\eta_{i,j}$ and a pole $\beta_0$.
Thus, $x(z_1, z_2)$ satisfying both equations of (\ref{2-1}) is given by
\begin{align}
x_i(z_1,z_2) &= (z_1-\alpha _0)^{-a_i} 
  \left( c_i+\sum^\infty_{j=1}d_{i,j}(\alpha _1,\cdots ,\alpha _{m-1}) (z_1-\alpha _0)^j\right) 
\label{4-10} \\
&= (z_1-\alpha _0)^{-a_i} \left( c_i+\sum^\infty_{j=1}d_{i,j}(y_1,\cdots ,y_{m-1})
\frac{(z_1-\alpha _0)^j}{(z_2-\beta_0)^{j}} \right),
\label{4-11}
\end{align}
where we used $d_{i,j}(\lambda \cdot A) = \lambda ^j d_{i,j}(A)$.
This gives a solution of (\ref{2-1}) as a function of $(z_1, z_2)$ as long as 
the right hand side converges.
Assume that the series (\ref{4-9}) converges when $|z_2 - \beta_0| \leq \varepsilon _2$
and (\ref{4-10}) converges when $|z_1 - \alpha _0| \leq \varepsilon_1$.
Let $\varepsilon $ be a small number.
Now we consider the solution restricted on the line
\begin{equation}
z_1- \alpha _0 = \varepsilon (z_2- \beta_0) \quad 
 \Leftrightarrow \quad  z_2 := q(z_1)= \frac{z_1 + \varepsilon \beta_0 -k_2}{\varepsilon +k_1},
\end{equation}
This yields
\begin{align}
x_i(z_1,q(z_1)) &= 
(z_1-\alpha _0)^{-a_i} \left( c_i+\sum^\infty_{j=1}d_{i,j}(y_1,\cdots ,y_{m-1}) \cdot \varepsilon ^j \right),
\label{4-12} \\
y_i &= \xi_i + \sum^\infty_{j=1} \eta_{i,j} (z_1-\alpha _0)^{j}/\varepsilon ^j. \nonumber
\end{align}
This is a convergent series when $|\varepsilon | < |\varepsilon _1|$
and $|z_1 - \alpha _0| < \varepsilon \varepsilon _2$.
Expanding it gives a new Laurent series solution
\begin{equation}
x_i(z_1,q(z_1)) = 
(z_1-\alpha _0)^{-a_i} \left( 
c_i+\sum^\infty_{j=1}d_{i,j}(\xi_1,\cdots ,\xi_{m-1})\cdot \varepsilon ^j + O(z_1-\alpha _0) \right),
\label{4-13}
\end{equation}
with the indicial locus
\begin{align}
(c_1 + \sum^\infty_{j=1}d_{1,j}(\xi_1, \cdots ,\xi_{m-1})\varepsilon ^j\,, \cdots ,\,
c_m + \sum^\infty_{j=1}d_{m,j}(\xi_1, \cdots ,\xi_{m-1})\varepsilon ^j).
\end{align}
If there are $k-2$ nonnegative integers among $\rho_2,\cdots ,\rho_{m-1}$,
(\ref{4-13}) contains $k-1$ free parameters 
($\alpha _0$ and $k-2$ parameters in $\eta_{i,j}$).
Suppose $\xi$ is a principal indicial locus of (\ref{4-8});
all $\rho_2, \cdots , \rho_{m-1}$ are nonnegative integers, then (\ref{4-13}) represents 
$m-1$ parameter family of Laurent series.
This is not a solution of $dx/dz_1 = F$ but a combination of $F$ and $G$;
\begin{align}
\frac{d}{dz_1}x_i(z_1, q(z_1)) 
  &= \frac{\partial }{\partial z_1}\Bigl|_{z_2 = q(z_1)}x_i(z_1, z_2) 
    + \frac{\partial }{\partial z_2}\Bigl|_{z_2 = q(z_1)} x_i(z_1, z_2)\cdot \frac{dq}{dz_1} \nonumber \\
&= f_i(x(z_1, q(z_1))) + \frac{1}{\varepsilon +k_1}g_i(x(z_1, q(z_1))).
\label{4-16}
\end{align}
This implies that there exists a lower indicial locus of the vector field $F+G/(\varepsilon +k_1)$ whose
K-exponents are given by $\rho_0=-1,\rho_1 = -1$ and $\rho_2, \cdots ,\rho_{m-1}$.
Since the number of free parameters is smaller than $m$, $x(z_1, q(z_1)) $ is a \textit{non}-principal
Laurent series solution.
\vskip.5\baselineskip

\textbf{Proposition \thedef.}
Suppose $\varepsilon $ and $|z_1-\alpha _0|$ are sufficiently small.
The Laurent series solution (\ref{4-13}) is convergent and it satisfies (\ref{4-16}).
In particular, there exists an indicial locus of the vector field $F+G/(\varepsilon +k_1)$ whose
K-exponents are given by $\rho_0=-1,\rho_1 = -1$ and $\rho_2, \cdots ,\rho_{m-1}$ for any small $\varepsilon \neq 0$.
\vskip.5\baselineskip

Although the series (\ref{4-13}) does not converge when $\varepsilon $ is large,
K-exponents $\rho_0=-1,\rho_1 = -1,\, \rho_2, \cdots ,\rho_{m-1}$ can be defined 
as the eigenvalues of the K-matrix of $F+G/(\varepsilon +k_1)$.
Since they analytically depend on $\varepsilon $ and are constants in $\varepsilon $
when $\varepsilon $ is small, we can take $\varepsilon \to \infty$
and obtain the main theorem in this subsection.
\vskip.5\baselineskip

\textbf{Theorem \thedef.}
There exists a lower indicial locus of the vector field $F$ whose
K-exponents are given by $\rho_0=-1,\rho_1 = -1$ and $\rho_2, \cdots ,\rho_{m-1}$.
\vskip.5\baselineskip

The expression of the corresponding indicial locus 
$c_i+\sum^\infty_{j=1}d_{i,j}(\xi_1,\cdots ,\xi_{m-1})\varepsilon ^j$
makes sense by an analytic continuation when $\varepsilon $ is large.
It is notable that $F$ has a lower indicial locus, but $G$ need not have.
When $\gamma =1$, the difference of assumptions for $F$ and $G$ are only (A3).
This suggests that (A3) is essential for the existence of lower indicial loci.
This is illustrated in the next example.
\vskip.5\baselineskip

\textbf{Example \thedef.}
Let us consider the two Hamiltonian functions
\begin{equation}
\left\{ \begin{array}{ll}
H_F(q_1, p_1, q_2, p_2) = (p_1^2/2-2 q_1^3)+(p_2^2/2-2 q_2^3),  &  \\
H_G(q_1, p_1, q_2, p_2) =  p_1^2/2-2 q_1^3,  \\
\end{array} \right.
\end{equation}
and let $F$ and $G$ be the corresponding Hamiltonian vector fields.
The vector field $F$ is a direct product of the two-dimensional system given in Example 2.8.
The weight is $(a_1, b_1, a_2, b_2) = (2,3,2,3)$ and $\gamma =1$.
The only $F$ satisfies the condition (A3).
The vector field $F$ has three indicial loci with K-exponents as 
\begin{align*}
(\text{P}_1) &: (q_1, p_1, q_2, p_2) =(1,-2,0,0), \quad \kappa = -1,2,3,6, \\
(\text{P}_2) &: (q_1, p_1, q_2, p_2)= (0,0,1,-2), \quad \kappa = -1,2,3,6,\\
(\text{P}_3) &:(q_1, p_1, q_2, p_2)= (1,-2,1,-2), \quad \kappa = -1,-1,6,6.
\end{align*}
The only $(\text{P}_3)$ is a lower locus.
For the first one $(\text{P}_1)$, the Laurent series solution of $(q_1, p_1)$ is given by (\ref{ex1}),
which has a free parameter $d_{1,6}$ at order $6$ (counting from $T^{-3}$).
$(q_2, p_2)$ is a holomorphic solution.
Thus, free parameters are the constant terms of the Taylor expansion 
those are at order 2 and 3, counting from $T^{-2}$ and $T^{-3}$, respectively.
These three parameters correspond to $\kappa = 2,3,6$.
Similarly for $(\text{P}_2)$.
For $(\text{P}_3)$, both of $(q_1, p_1)$ and $(q_2, p_2)$ are Laurent series solutions
of the form (\ref{ex1}). Hence, there are two free parameters at order $6$,
that correspond to $\kappa = 6,6$.

For $(\text{P}_1)$, the flow of the free parameters are given by
\begin{align*}
\alpha _0' = -1,\,\, \alpha _1' = -\alpha _2,\,\, \alpha _2' = -6 \alpha _1^2,\,\, \alpha _3' = 0.
\end{align*}
It has only one indicial locus and its K-exponents are given by
$\kappa = -1, -1, 6, 6$, that confirms Theorem 4.6.

On the other hand, the vector field $G$ has only one indicial locus whose K-exponents
are given by $(-1,2,3,6)$.
There are no lower indicial loci.
When $\gamma =1$, we cannot distinguish $F$ and $G$ by assumptions (A1) and (A2) in Section 2.
This suggests that the assumption (A3) ``$F(x) =0$ only when $x=0$"
plays a crucial role for the existence of a lower indicial locus,
though the authors do not know a clear statement.

\begin{figure}
\begin{center}
\includegraphics[scale=0.6]{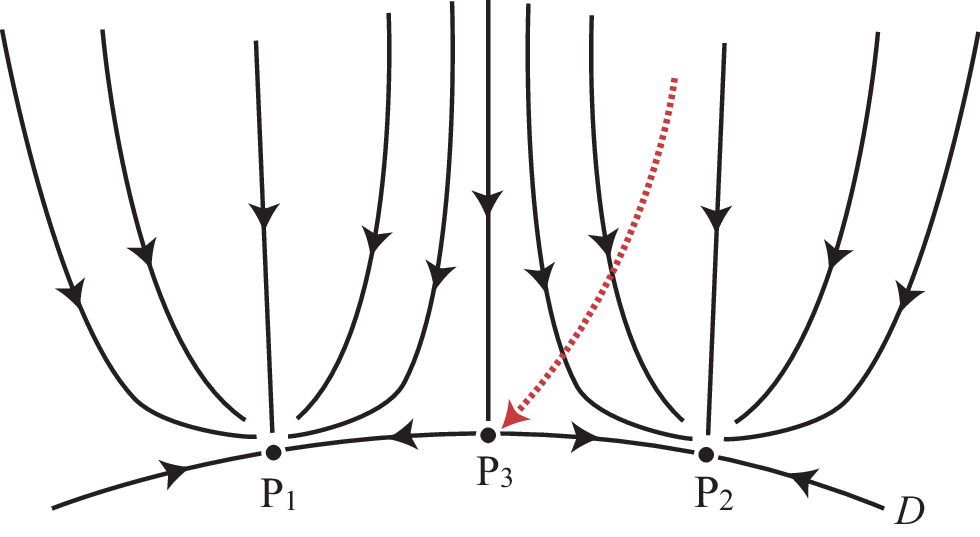}
\caption[]{A schematic view of the flow of $F$ on the compactified space.
The original phase space $\C^4$ is compactified by attaching $D$ at ``infinity".
There are three fixed points on $D$ that correspond to three indicial loci.
The red dotted orbit indicates an orbit of the vector field $F+G/(\varepsilon +k_1)$ given in Proposition 4.5.
}
\label{fig1}
\end{center}
\end{figure}

Fig. \ref{fig1} represents a schematic view of the flow of $F$.
In Chiba \cite{Chi1, Chi4}, the geometry of families of Laurent series solutions
of quasi-homogeneous vector fields are investigated via the weighted projective space (in this example,
it is $\C P^4(2,3,2,3)$).
This is a compact manifold constructed by attaching an $m-1$ dimensional manifold, 
denoted by $D$ in the figure, to the original phase space $\C^m$ at infinity.
A given vector field on $\C^m$ is extended to a vector field on $\C^m \cup D$,
and it is shown that there is a one-to-one correspondence between indicial loci and 
fixed points on $D$.
Further, eigenvalues of the Jacobi matrix of $F$ at the fixed point on $D$
coincide with K-exponents.
In the figure, $(\text{P}_1)$ and $(\text{P}_2)$ are stable fixed points that 
correspond to principal indicial loci, respectively, and $(\text{P}_3)$ is a saddle point 
that corresponds to the lower indicial locus.
The space is divided into two regions that are occupied by two families of 
principal Laurent series solutions, 
and their boundary is occupied by the non-principal Laurent series solutions.
The red dotted orbit indicates an orbit of the vector field $F+G/(\varepsilon +k_1)$
given in Proposition 4.5.

\subsection{$\gamma \geq 2$}

Let us consider the case $\gamma \geq 2$.
For the equation $d\alpha _0/dz_2 = \hat{g}_0(A)$,\, the right hand side is not a constant
because of Lemma 4.2, however, it may become identically zero.
For it, we can prove the next proposition.
The proof is given in Appendix.
\vskip.5\baselineskip

\textbf{Proposition \thedef.}
Let $c$ be an isolated principal indicial locus of $F=(f_1,\cdots ,f_m)$ and 
$\kappa_1 \geq 1$ be the smallest K-exponent other than $\kappa_0 = -1$.
If the vector $(d_{1,\kappa_1},\cdots ,d_{m, \kappa_1})$ is not an eigenvector associated with a
zero eigenvalue of the Jacobi matrix of $G$ at $c$, then $\hat{g}_0(A)$ is not identically zero.
\vskip.5\baselineskip

In what follows, we assume that $\hat{g}_0(A) \not\equiv 0$, not identically zero.
We rewrite (\ref{4-6}) as
\begin{equation}
\frac{d\alpha _i}{d\alpha _0} = \frac{\hat{g}_i(A)}{\hat{g}_0(A)}, \quad i=1,\cdots ,m-1,
\quad A=(\alpha _1,\cdots ,\alpha _{m-1}).
\label{4-20}
\end{equation}
Because of (\ref{4-5}), it has the degree $1$ with respect to the weight $(\kappa _1,\cdots ,\kappa _{m-1})$.
\vskip.5\baselineskip

\textbf{Lemma \thedef.} Let $\xi = (\xi_0,\cdots ,\xi_{m-1})$ be an indicial locus of (\ref{4-6}).
Assume that $\hat{g}_0(\xi) \neq 0$.
Then, (\ref{4-20}) has an indicial locus $\widetilde{\xi} = (\widetilde{\xi}_1 ,\cdots , \widetilde{\xi}_{m-1})$ 
with $\widetilde{\xi}_i = \xi_0^{\kappa _i} \xi_i$.
\vskip.5\baselineskip

\textbf{Proof.}
Since (\ref{4-6}) has the degree $\gamma $ with respect to the weight $\kappa_0, \cdots , \kappa_{m-1}$,
the indicial locus is given by the root of the equation (see (\ref{3-10}))
\begin{equation*}
\hat{g}_0(\xi) = -\kappa _0 \xi_0/\gamma = \xi_0 / \gamma ,\quad \hat{g}_i(\xi) = -\kappa _i \xi_i /\gamma.
\end{equation*}
Putting $\xi_i = \xi_0^{-\kappa _i}\widetilde{\xi}_i$ gives
\begin{align}
\xi_0 \frac{\hat{g}_i (\xi)}{\hat{g}_0 (\xi)} = -\kappa_i \xi_i
\quad & \Rightarrow \quad 
  \xi_0\frac{\xi_0^{-\kappa_i-\gamma }
  \hat{g}_0(\widetilde{\xi})}{\xi_0^{-\kappa_0 - \gamma }\hat{g}_0 (\widetilde{\xi})} 
  = - \kappa_i \xi_0^{-\kappa_i}\widetilde{\xi}_i \nonumber \\
& \Rightarrow \quad 
  \frac{\hat{g}_i(\widetilde{\xi})}{\hat{g}_0 (\widetilde{\xi})} = - \kappa_i \widetilde{\xi}_i
\label{4-21}
\end{align}
Hence, $\widetilde{\xi}$ satisfies the definition of an indicial locus of (\ref{4-20}). $\hfill \Box$
\vskip.5\baselineskip

\textbf{Proposition \thedef.}
Let $\rho_0 = -1,\, \rho_1 = -1/\gamma , \rho_2,\cdots ,\rho_{m-1}$ be K-exponents of (\ref{4-6}) at a locus $\xi$
(see Lemma 4.3).
Then, the K-exponents of (\ref{4-20}) at the locus $\widetilde{\xi}$ above are given by $\rho_0 =-1$
and $\gamma \rho_2 ,\cdots ,\gamma \rho_{m-1}$.
\vskip.5\baselineskip

\textbf{Proof.} Denote
\begin{align*}
\mathrm{grad}(\hat{g}_0) =  \left( \frac{\partial \hat{g}_0}{\partial \alpha _1}(\xi) ,\cdots , 
     \frac{\partial \hat{g}_0}{\partial \alpha _{m-1}}(\xi)\right),\quad
(J\hat{g}) = \left\{ \frac{\partial \hat{g}_i}{\partial \alpha _j} (\xi)\right\}^{m-1}_{i,j=1},
\end{align*}
and $\mathrm{diag}(\kappa ) = \mathrm{diag} (\kappa _1,\cdots ,\kappa _{m-1})$.
Then, the K-matrix of (\ref{4-6}) given in (\ref{4-7}) is written as
\begin{equation}
K_\gamma (\xi) = \left(
\begin{array}{@{\,}c|c@{\,}}
0 & \mathrm{grad}(\hat{g}_0) \\ \hline
0& (J\hat{g})
\end{array}
\right) + \left(
\begin{array}{@{\,}c|c@{\,}}
-1/\gamma & 0 \\ \hline
0 & \mathrm{diag}(\kappa )/\gamma 
\end{array}
\right).
\end{equation}
Further, put $a = 1/(\gamma  \hat{g}_0(\xi))$ and $v = (\kappa _1 \xi_1 ,\cdots , \kappa _{m-1}\xi_{m-1})^T$,
which is an eigenvector associated with the K-exponent $\rho_0 = -1$ of the subsystem (\ref{4-8}).
Define
\begin{equation}
P = \left(
\begin{array}{@{\,}c|c@{\,}}
a & 0 \\ \hline
-v & I 
\end{array}
\right).
\end{equation}
Then, we can verify that 
\begin{equation*}
P^{-1} K_\gamma (\xi)P = \left(
\begin{array}{@{\,}c|c@{\,}}
k_{11} & \mathrm{grad}(\hat{g}_0) / a \\ \hline
k_{21} & k_{22} 
\end{array}
\right),
\end{equation*}
where 
\begin{equation*}
\left\{ \begin{array}{l}
k_{11} = -\mathrm{grad}(\hat{g}_0)\cdot v/a - 1/\gamma ,  \\
k_{21} = -(\mathrm{grad} (\hat{g}_0) \cdot v)v / a- v/\gamma - (J\hat{g})v - \mathrm{diag}(\kappa ) v/\gamma ,  \\
k_{22} = (J\hat{g}) + v\cdot \mathrm{grad}(\hat{g}_0)/a + \mathrm{diag}(\kappa )/\gamma .
\end{array} \right.
\end{equation*}
Since $v$ is an eigenvector of the K-exponent $-1$ of (\ref{4-8}), we have
\begin{align*}
((J\hat{g}_0) + \mathrm{diag}(\kappa )/\gamma )v = -v.
\end{align*}
From the derivative of (\ref{4-5}) for $l=0$ at $\lambda =1$, we have
\begin{align*}
\sum^{m-1}_{j=1} \frac{\partial \hat{g}_0}{\partial \alpha _j}(\xi) \kappa _j \xi_j = (-1+\gamma ) \hat{g}_0(\xi).
\end{align*}
By using these two equalities, we can show that $k_{11} = -1$ and $k_{21} = 0$.
Thus, the eigenvalues of $k_{22}$ are $\rho_1 = -1/\gamma $ and $\rho_2, \cdots , \rho_{m-1}$.
This implies that the eigenvalues of $\gamma \cdot k_{22}$ are $-1$ and $\gamma \rho_2, \cdots , \gamma \rho_{m-1}$.
Hence, it is sufficient to prove that the K-matrix of (\ref{4-20}) at $\widetilde{\xi}$ is conjugate to $\gamma \cdot k_{22}$.

By a straightforward calculation with the aid of (\ref{4-21}) and (\ref{3-3}), it is easy to see that
the $(i,j)$-component of the K-matrix of (\ref{4-20}) is given by
\begin{align*}
\widetilde{K}(\widetilde{\xi})_{i,j}
 & = \frac{1}{\hat{g}_0 (\widetilde{\xi})} 
   \left( \frac{\partial \hat{g}_i}{\partial \alpha _j}(\widetilde{\xi})
 + \kappa _i \widetilde{\xi}_i \frac{\partial \hat{g}_0}{\partial \alpha _j}(\widetilde{\xi})  \right) 
  + \kappa _i\delta _{ij} \\
&= \frac{\xi_0 }{\hat{g}_0 (\xi)} 
   \left( \xi_0 ^{\kappa _i-\kappa _j}\frac{\partial \hat{g}_i}{\partial \alpha _j}(\xi)
 + \xi_0 ^{\kappa _i-\kappa _j-1} \kappa _i \xi_i \frac{\partial \hat{g}_0}{\partial \alpha _j}(\xi)  \right) 
  + \kappa _i\delta _{ij}.
\end{align*}
Put $Q = \mathrm{diag} (\xi_0 ^{\kappa _1} , \cdots , \xi_0 ^{\kappa _{m-1}})$.
Then, we can show that $Q^{-1}\widetilde{K}(\widetilde{\xi}) Q = \gamma k_{22}$. $\hfill \Box$
\vskip.5\baselineskip

Let $\xi = (\xi_0,\cdots ,\xi_{m-1})$ be an indicial locus of (\ref{4-6}) and 
assume that the associated Puiseux series solution contains $m-1$ free parameters.
As is shown in the end of Sec.3, its K-exponents should be $\rho_i \in \mathbb{N}/\gamma $
for $i=2,\cdots ,m-1$.
Then, $\widetilde{\xi}$ is a principal indicial locus of the system (\ref{4-20})
with K-exponents $-1, \gamma \rho_2, \cdots ,\gamma \rho_{m-1}$
satisfying $\gamma \rho_2, \cdots ,\gamma \rho_{m-1} \in \N$.
Since (\ref{4-20}) is a rational vector field of degree $1$ with respect to the weight $(\kappa _1,\cdots ,\kappa _{m-1})$,
it has a Laurent series solution of the form
\begin{align}
\alpha _i(\alpha _0) = (\alpha _0-\beta_0)^{-\kappa_i} 
  \left( \widetilde{\xi}_i+\sum^\infty_{j=1} \widetilde{\eta}_{i,j} (\alpha _0-\beta_0) ^j\right)
 =: (\alpha _0-\beta_0)^{-\kappa_i}\widetilde{y}_i,
\label{4-22}
\end{align}
for $i=1,\cdots ,m-1$.
Free parameters are contained in $\widetilde{\eta}_{i,j}$ for $j= \gamma \rho_2, \cdots ,\gamma \rho_{m-1}$.
Hence, $x(z_1, \alpha _0)$ satisfying both equations of (\ref{2-1}) is given by
\begin{align}
x_i(z_1, \alpha _0) &= (z_1-\alpha _0)^{-a_i} 
  \left( c_i+\sum^\infty_{j=1}d_{i,j}(\alpha _1,\cdots ,\alpha _{m-1}) (z_1-\alpha _0)^j\right) 
\label{4-23} \\
&= (z_1-\alpha _0)^{-a_i} \left( c_i+\sum^\infty_{j=1}d_{i,j}(\widetilde{y}_1,\cdots ,\widetilde{y}_{m-1})
\frac{(z_1-\alpha _0)^j}{(\alpha _0-\beta_0)^{j}} \right),
\end{align}
where $\alpha _0 = \alpha _0(z_2)$ is related to $z_2$ though $d\alpha _0/dz_2 = \hat{g}_0(A)$.
Assume that the series (\ref{4-22}) converges when $|\alpha _0 - \beta_0| \leq \varepsilon _2$
and (\ref{4-23}) converges when $|z_1 - \alpha _0| \leq \varepsilon_1$.
Let $\varepsilon $ be a small number.
Now we consider the solution restricted on the line
\begin{equation}
z_1- \alpha _0 = \varepsilon (\alpha_0 - \beta_0) \quad 
 \Leftrightarrow \quad  \alpha _0 := q(z_1)= \frac{z_1 + \varepsilon \beta_0}{1+\varepsilon},
\end{equation}
This yields
\begin{align}
x_i(z_1,q(z_1)) &= 
(z_1-\alpha _0)^{-a_i} 
  \left( c_i+\sum^\infty_{j=1}d_{i,j}(\widetilde{y}_1,\cdots ,\widetilde{y}_{m-1}) \cdot \varepsilon ^j \right),
\label{4-26} \\
\widetilde{y}_i &= \widetilde{\xi}_i + \sum^\infty_{j=1} \widetilde{\eta}_{i,j} (z_1-\alpha _0)^{j}/\varepsilon ^j. \nonumber
\end{align}
This is a convergent series when $|\varepsilon |<|\varepsilon _1|$ and $|z_1 - \alpha _0| < \varepsilon \varepsilon _2$.
Expanding it gives a new Laurent series solution as in Section 4.1.
If all $\gamma \rho_2, \cdots , \gamma \rho_{m-1}$ are positive integers, then it represents 
$m-1$ parameter family of Laurent series. This satisfies
\begin{align*}
\frac{d}{dz_1}x_i(z_1, q(z_1)) 
  &= \frac{\partial }{\partial z_1}\Bigl|_{\alpha _0= q(z_1)}x_i(z_1, \alpha _0) 
    + \frac{\partial }{\partial \alpha _0}\Bigl|_{\alpha _0 = q(z_1)} x_i(z_1, \alpha _0) \cdot \frac{dq}{dz_1} \nonumber \\
&= f_i(x(z_1, q(z_1))) 
+ \frac{1}{1+\varepsilon}\frac{\partial }{\partial \alpha _0}\Bigl|_{\alpha _0 = q(z_1)} x_i(z_1, \alpha _0).
\end{align*}
On the other hand, we have
\begin{align*}
g_i(x)= \frac{\partial x_i}{\partial z_2} = \frac{\partial x_i}{\partial \alpha _0}\cdot \frac{d\alpha _0}{dz_2}
= \frac{\partial x_i}{\partial \alpha _0}\cdot \hat{g}_0(A).
\end{align*}
This shows
\begin{equation}
\frac{d}{dz_1}x_i(z_1, q(z_1)) = 
f_i(x(z_1, q(z_1))) 
+ \frac{1}{1+\varepsilon} \frac{g_i(x(z_1, q(z_1)))}{\hat{g}_0(A)}.
\label{4-27}
\end{equation}
Here, $\hat{g}_0(A)$ is regarded as a function of $x = x(z_1, q(z_1))$ through $A = \Phi^{-1}(x)$.
By the same way as Section 4.1, we obtain the next results.
\vskip.5\baselineskip

\textbf{Proposition \thedef.}
Suppose $\varepsilon $ and $|z_1-\alpha _0|$ are sufficiently small.
The Laurent series solution (\ref{4-26}) is convergent and it satisfies (\ref{4-27}).
In particular, there exists an indicial locus of the vector field $F+G/(\hat{g}_0\cdot (1+\varepsilon))$ whose
K-exponents are given by $\rho_0=-1,\rho_1 = -\gamma $ 
and $\gamma \rho_2, \cdots ,\gamma \rho_{m-1}$ for any small $\varepsilon \neq 0$.
\vskip.5\baselineskip

\textbf{Theorem \thedef.}
There exists a lower indicial locus of the vector field $F$ whose
K-exponents are given by $\rho_0=-1,\rho_1 = -\gamma $ and $\gamma \rho_2, \cdots ,\gamma \rho_{m-1}$ .
\vskip.5\baselineskip
\vskip.5\baselineskip


\textbf{Example \thedef.}
Let us consider the two Hamiltonian functions
\begin{equation}
\left\{ \begin{array}{ll}
H_F(q_1, p_1, q_2, p_2) = 2p_1p_2 + 3p_2^2 q_1+q_1^4 - q_1^2q_2-q_2^2,  &  \\
H_G(q_1, p_1, q_2, p_2) = p_1^2 + 2p_1 p_2 q_1 - q_1^5 + p_2^2 q_2 + 3 q_1^3 q_2-2 q_1 q_2^2,  \\
\end{array} \right.
\end{equation}
and let $F$ and $G$ be the corresponding Hamiltonian vector fields.
This is known as the autonomous version of the 4-dimensional first Painlev\'{e} equation \cite{Chi4}.

The weight is $(a_1, b_1, a_2, b_2) = (2,5,4,3)$ and $\gamma =3$.
With this weight, the weighted degrees of Hamiltonian are $\deg(H_F) = 8$ and $\deg(H_G) = 10$.

The vector field $F$ has two indicial loci with K-exponents as 
\begin{align*}
(\text{P}_1) &: (q_1, p_1, q_2, p_2) =(1,1,1,-1), \quad \kappa = -1,2,5,8, \\
(\text{P}_2) &: (q_1, p_1, q_2, p_2)= (3,27,0,-3), \quad \kappa = -3,-1,8,10.
\end{align*}
The only $(\text{P}_2)$ is a lower locus.
It is known that for Hamiltonian vector fields, the K-exponents always appear as a pair
in the sense that;
\vskip.5\baselineskip

\textbf{Proposition \thedef \,\cite{Chi4}.}
For a quasi-homogeneous Hamiltonian system $F$ satisfying (\ref{A1}),
if $\kappa$ is a Kovalevskaya exponent, so is $\mu$ given by
 $\kappa + \mu = \deg(H_F)-1$.
Further, the following formula holds
\begin{equation}
\kappa + \mu = \deg(H_F)-1 = \deg (q_i) + \deg (p_i).
\end{equation}
\vskip.5\baselineskip

In this example, it means that
\begin{align*}
\kappa_{i} + \kappa_{4-i-1} = \deg(q_j) + \deg(p_j) = \deg(H_F)-1 = 7,\quad
i=0,1,\,\, j=1,2.
\end{align*}
Applying the proposition to $\kappa_0 = -1$, it turns out that $\deg(H_F)$
is always a K-exponent.
By Theorem 4.12, $\kappa_1 = -\gamma $ is a K-exponent for a lower indicial locus.
In this example, this means $(-3)+\kappa_2 = 7$ and $\kappa_2 = 10 = \deg (H_G)$.
Thus, we can obtain K-exponents of lower indicial locus from the weight
of Hamiltonian functions.

For $(\text{P}_1)$, the flow of the free parameters are given by
\begin{align*}
\alpha _0' = 3 \alpha _1,\,\, \alpha _1' = -\frac{3}{2}\alpha _2,\,\, 
  \alpha _2' = -54 \alpha _1^4,\,\, \alpha _3' = 42 \alpha _1^3 \alpha _2.
\end{align*}
It has three lower indicial loci and their K-exponents are given by
$\rho = -1/3, -1, 8/3, 10/3$, that confirms Theorem 4.12.
The vector field $G$ has a lower indicial locus whose K-exponents are also
$\rho = -1/3, -1, 8/3, 10/3$. 
As was explained in Remark 4.4, it is induced from the lower locus of the parameter flow.

\appendix
\section{The proof of Proposition 4.8.}

We again state the proposition;
\vskip.5\baselineskip

\textbf{Proposition \thedef.}
Let $c$ be an isolated principal indicial locus of $F=(f_1,\cdots ,f_m)$ and 
$\kappa_1 \geq 1$ be the smallest K-exponent other than $\kappa_0 = -1$.
If the vector $(d_{1,\kappa_1},\cdots ,d_{m, \kappa_1})$ is not an eigenvector associated with a
zero eigenvalue of the Jacobi matrix of $G$ at $c$, then $\hat{g}_0(A)$ is not identically zero.
\vskip.5\baselineskip

\textbf{Proof.}
The first half of the proof is similar to that of Prop.4.1.
We repeat it to fix our notation.
A Laurent series solution of $dx_i/dz_1 = f_i(x)$ is given by
\begin{align*}
x_i &= (z_1 - \alpha _0)^{-a_i} 
  \left( c_i + \sum_{j=1} d_{i,j}(A) (z_1 - \alpha _0) ^j\right), \\
&=:  (z_1 - \alpha _0)^{-a_i} y_i(A), \quad A=(\alpha _1(z_2), \cdots ,\alpha _{m-1}(z_2)) .
\end{align*}
Substituting it into $dx_i/dz_2 = g_i(x)$, for the left hand side we have
\begin{align*}
\frac{dx_i}{dz_2} &= a_ic_i(z_1-\alpha _0)^{-a_i-1} \frac{d\alpha _0}{dz_2}
  + \sum_{j=1}(a_i-j)d_{i,j}(A) (z_1-\alpha _0)^{j-a_i-1} \frac{d\alpha _0}{dz_2} \\
  & \qquad\qquad + \sum_{j=1}\frac{d}{dz_2}\left( d_{i,j}(A)\right) \cdot (z_1-\alpha _0)^{j-a_i}.
\end{align*}
For the right hand side,
\begin{align*}
g_i(x) &= (z_1 - \alpha _0)^{-a_i-\gamma } g_i(y_1, \cdots , y_m) \\
&=:  (z_1 - \alpha _0)^{-a_i-\gamma } \sum_{k=0} g_{i,k}(A)(z_1-\alpha _0)^k,
\end{align*}
where $g_{i,k}(A)$ is a coefficient of the Taylor expansion of $g_i(y_1, \cdots , y_m)$.
They are given through
\begin{align}
g_i(y_1, \cdots , y_m)
  = g_i(c) + \sum^m_{l=1} \frac{\partial g_i}{\partial x_l}(c) \cdot
    \left( \sum_{j=1} d_{l,j}(A) (z_1-\alpha_0)^j \right) + \cdots .
\label{A-1}
\end{align}
Comparing coefficients of both sides of $dx_i/dz_2 = g_i(x)$, we obtain
$g_{i,0} = \cdots = g_{i,\gamma -2} = 0$ and
\begin{align*}
(z_1-\alpha _0)^{-a_i-1} &:\quad  a_ic_i \frac{d\alpha _0}{dz_2} = g_{i, \gamma -1}(A) \\
(z_1-\alpha _0)^{-a_i} &:\quad (a_i-1) d_{i,1}(A)\frac{d\alpha _0}{dz_2} = g_{i,\gamma }(A) \\
(z_1-\alpha _0)^{-a_i+1} &:\quad (a_i-2) d_{i,2}(A)\frac{d\alpha _0}{dz_2}
    + \frac{d}{dz_2}d_{i,1}(A)  = g_{i,\gamma +1 }(A)\\
&\quad\quad\quad \vdots \\
(z_1-\alpha _0)^{-a_i+j} &:\quad (a_i-j-1) d_{i,j+1}(A)\frac{d\alpha _0}{dz_2}
    + \frac{d}{dz_2}d_{i,j}(A)  = g_{i,\gamma +j }(A) \\
&\quad\quad\quad \vdots 
\end{align*}
Now we assume that $\hat{g}_0(A) = g_{i, \gamma -1}(A)/(a_ic_i) = 0$ and we will derive a contradiction.
Since $d\alpha _0/dz_2 = \hat{g}_0(A)= 0$, the above equations yield $g_{i,\gamma }(A) = 0$ and
\begin{align}
\frac{d}{dz_2}d_{i,j}(A)  = g_{i,\gamma +j }(A) , \quad j=1,2,\cdots .
\label{A-2}
\end{align}
Let $\kappa_1\geq 1$ be the smallest K-exponent.
Then, $d_{i,j}(A) = 0$ for $j=1,\cdots ,\kappa_1-1$ and $i=1,\cdots ,m$, and 
$d_{l,\kappa_1}(A) = \alpha _1$ for some $l$ by the definition of the free parameter $\alpha _1$
(see Definition 2.4 below).
Hence, (\ref{A-1}) becomes
\begin{align*}
g_i(y_1, \cdots ,y_m) = g_i(c) + \sum^m_{l=1}\frac{\partial g_i}{\partial x_l}(c) \cdot 
\left( d_{l,\kappa_1} (A) (z_1-\alpha _0)^{\kappa_1} + \cdots \right) + O((z_1-\alpha _0)^{2\kappa_1}).
\end{align*}
This provides
\begin{align*}
g_{i,\kappa_1}(A) = \sum^m_{l=1}\frac{\partial g_i}{\partial x_l}(c)d_{l,\kappa_1} (A).
\end{align*}
Further, $d_{i,j}(A) = 0$ for $j=1,\cdots ,\kappa_1-1$ with (\ref{A-2}) gives
$g_{i, \gamma +1} = \cdots  = g_{i,\gamma + \kappa_1 -1} = 0$.
In particular, 
\begin{align*}
g_{i,\kappa_1}(A) = 0 = \sum^m_{l=1}\frac{\partial g_i}{\partial x_l}(c)d_{l,\kappa_1} (A).
\end{align*}
It contradicts with the assumption of the proposition.
$\hfill \Box$



\begin{thebibliography}{99}
\setlength{\baselineskip}{0pt}

\bibitem{Adl}
M. Adler, P. Moerbeke, P. Vanhaecke,
Algebraic Integrability, Painlev\'{e} Geometry and Lie Algebras,
Springer-Verlag, Berlin, (2004).

\bibitem{Bor}
A. V. Borisov, S. L. Dudoladov,,
Kovalevskaya exponents and Poisson structures,
Regul. Chaotic Dyn. 4 (1999), no. 3, 13-20. 

\bibitem{Chi1}
H. Chiba, 
Kovalevskaya exponents and the space of initial conditions of a quasi-homogeneous vector field,
J. Diff. Equ. 259, no. 12, 7681-7716, (2015).

\bibitem{Chi2}
H. Chiba,
The first, second and fourth Painlev\'{e} equations on weighted projective spaces,
 J. Diff. Equ. 260, no. 2, 1263-1313, (2015).

\bibitem{Chi3}
H. Chiba,
The third, fifth and sixth Painlev\'{e} equations on weighted projective spaces,
 SIGMA 12, 019 (2016).

\bibitem{Chi4}
H. Chiba,
Weights, Kovalevskaya exponents and the Painleve property,
Annales de l'Institut Fourier, 74, no. 2, 811-848, (2024).

\bibitem{Cho}
S. N. Chow, C. Li, D. Wang,
Normal forms and bifurcation of planar vector fields,
Cambridge University Press, Cambridge, (1994).

\bibitem{CosCos}
O. Costin, R. D. Costin,
Singular normal form for the Painlev\'{e} equation P1,
Nonlinearity 11 (1998),  no.5, 1195-1208.

\bibitem{Gor}
A. Goriely,
Painlev\'{e} analysis and normal forms theory,
Phys. D 152/153 (2001), 124-144.

\bibitem{Gor2}
A. Goriely,
Integrability and nonintegrability of dynamical systems,
Advanced Series in Nonlinear Dynamics, 19,
World Scientific Publishing Co., Inc., River Edge, NJ (2001).

\bibitem{HuYa}
J. Hu, M. Yan,
Painlev\'{e} Test and the Resolution of Singularities for Integrable Equations,
(arXiv:1304.7982).

\end{thebibliography}
\end{document}